\newif
\ifconfver
\confvertrue        
\ifconfver
    \documentclass[10pt,journal,final,twoside]{IEEEtran}
\else
    \documentclass[11pt,draftcls,onecolumn]{IEEEtran}
\fi
\usepackage{xspace,empheq,fancybox,amsmath,amssymb,graphicx,epstopdf,epsfig,subfigure,syntonly,times,amsthm} \usepackage{psfrag,color,bm,array,cite}
\usepackage{url,cite,footnote,xspace,syntonly,algorithm,algorithmic}
\usepackage{verbatim,multirow}
\usepackage[T1]{fontenc}
\usepackage{mathtools}
\usepackage{amsmath,amsfonts,amsthm,bm} 
\usepackage{graphicx}
\usepackage{mwe}
\usepackage{caption}
\usepackage{cite}
\usepackage{amsmath,amssymb,amsfonts}
\usepackage{algorithmic}
\usepackage{graphicx}
\usepackage{textcomp}
\usepackage{xcolor}
\usepackage{multicol, blindtext}
\usepackage{xspace,empheq,fancybox,amsmath,amssymb,graphicx,epstopdf,epsfig,syntonly,times,amsthm} \usepackage{psfrag,color,bm,array,cite}
\usepackage{url,cite,footnote,xspace,syntonly,algorithm,algorithmic}
\usepackage{verbatim,multirow}
\usepackage[T1]{fontenc}
\usepackage{mathtools}
\usepackage{subfig}
\usepackage{amsmath,amsfonts,amsthm,bm}
\usepackage{graphicx}
\usepackage{mwe}
\usepackage{caption}
\usepackage{stfloats}
\usepackage{amsfonts}
\usepackage{stackengine}
\usepackage{lipsum}

\input{mysymbol.sty}

\newcommand{\yuqi}{\color{black}{}}

\definecolor{orange}{rgb}{1,0.5,0}

\allowdisplaybreaks

\begin{document}


\title{Efficient Identification of Bus Split Events \\ Using Synchrophasor Data}
\author{
\IEEEauthorblockN{Yuqi Zhou},~\IEEEmembership{Student Member, IEEE}, and \IEEEauthorblockN{Hao Zhu},~\IEEEmembership{Senior Member, IEEE}

}

\renewcommand{\thepage}{}
\maketitle
\pagenumbering{arabic}

%
\begin{abstract}
Accurate grid topology information is of paramount importance for routine power system operations, while the growing availability of synchrophasor data offers the opportunity to identify topology changes in real time. Identification of bus split events, where the substation becomes electrically disconnected, is becoming increasingly important for maintaining the security of power systems. This paper aims to provide an efficient modeling and monitoring framework for bus split events {\yuqi by using a concise bus-branch representation}. The linear sensitivity analysis is first performed to quickly evaluate the grid-wide impact of such events. Furthermore, the synchrophasor data enabled identification problem is formulated by matching the changes in bus phase angles (and possibly line flows). To address the resultant bilinear multiplication involving the binary connectivity variables, the McCormick relaxation technique is leveraged to attain an equivalent mixed-integer linear program reformulation that is efficiently solvable. Numerical studies on the IEEE 14-bus and 300-bus systems demonstrate the validity and efficiency of the proposed identification algorithm towards real-time implementation.

\end{abstract}




\section{Introduction}\label{sec:intro}

\IEEEPARstart{A}{ccurate} information of grid topology is crucial for performing various power system operation and maintenance tasks \cite[Ch. 12]{wood2013power}. Recently, the power grids have witnessed more frequent occurrence of unintentional circuit breaker (CB) actions, due to either misoperations of substation protection systems under renewable generation \cite{EPRI01} or  malicious attacks as in the 2015/16 Ukrainian blackouts \cite{case2016analysis}. Changes of CB status not only lead to the disconnection of power lines or generation/loads, but also can  give rise to the \textit{bus split} events when a substation's bus bars become disconnected \cite{abur1995identifying}. Recent reports \cite{kassakian2011future,WECC1,NERC2,NERC1} point out the increasing importance of developing efficient grid topology modeling and monitoring techniques that can include bus split events.

{\yuqi
Traditionally, grid topology is processed and updated by the energy management system (EMS) using the Supervisory Control and Data Acquisition (SCADA) inputs \cite[Ch. 1]{abur2004power}. A compact bus-branch model is typically obtained by examining the statuses of CBs and switching devices. \textit{Topology errors} due to incorrect reports of CB statuses can be identified by the \textit{generalized} state estimator using SCADA data \cite{wu1989detection,clements1998topology,exposito2001reduced,korres2002identification,lourencco2004bayesian,lourencco2014topology}.
Recent deployment of phasor measurement units (PMUs) enables timely topology identification using the synchrophasor data of high resolution and sampling rate. 
Fast algorithms for identifying (even multiple) transmission line outages has been developed in \cite{tate2008line,zhu2012sparse,emami2012external,chen2015quickest}, and also for distribution grids \cite{deka2016estimating,cavraro2018graph,ardakanian2019identification}. {\yuqi Furthermore, model-free approaches using synchrophasor data to detect and identify general disturbance events have been considered recently. For example, various techniques have been advocated for this purpose, such as subspace analysis \cite{wu2016online}, ensemble learning \cite{zhou2018ensemble}, graph signal processing \cite{ramakrishna2019detection}, and wavelet transform \cite{wang2019data}. Nonetheless, none of these methods has specifically accounted for bus split events. } 

The problem of identifying CB status and bus split events is known to incur high complexity and identifiability issues. To capture  CB statuses and substation connectivity, existing approaches typically use the detailed node-breaker model which increases the problem size and computational time; see e.g.,   \cite{heidarifar2015network,park2019sparse,donmez2020parallel}. 
For example, a robust state estimation (SE) has been proposed in \cite{donmez2020parallel} to account for the SE topology errors, yet it needs to incorporate the CB flows as additional state variables.
Detecting CB actions from the synchrophasor data has been considered in \cite{kekatos2012joint}, where lack of identifiability has been observed therein due to the complicated node-breaker model. 
Thus, the problem of efficiently analyzing and identifying bus split events still remains open.   }



The goal of this paper is to {\yuqi develop an efficient modeling and monitoring framework for  bus split events by incorporating them into the concise bus-branch representation}. For the modeling part, we approach the dc power flow model \cite{stott2009dc} based sensitivity analysis through simplifying the matrix computations. The sensitivity analysis results allow for constructing an equivalent bus-branch model for the post-split system, and are essential for executing fast contingency screening and security-constrained economic dispatch tasks; see e.g., \cite[Ch. 11]{wood2013power}. 
{\yuqi Such concise bus-branch model for bus split events further allows for synchrophasor data based identification by searching for the post-split scenario of the best match to the measurements. Hence, one could enumerate all possible bus split scenarios and use the aforementioned sensitivity analysis to obtain the resultant post-split grid responses.} Nonetheless, the number of scenarios in this exhaustive search approach would grow exponentially with the number of substation connections, limiting it from real-time implementation. To address this complexity issue, we formulate the identification problem with binary variables indicating substation connections, for which the popular \textit{McCormick relaxation} technique \cite{McCormick1976computability} can be adopted to replace the \textit{nonconvex} bilinear terms with equivalent \textit{convex} linear inequalities. This reformulation  leads to a tractable mixed-integer linear programming problem that is efficiently solvable  for real-time identification. {\yuqi By developing a concise representation of  bus split events, our work directly addresses the lack of consideration and complexity issues in existing framework for grid analysis and monitoring}.

The rest of the paper is organized as follows. Section \ref{sec:ps} introduces the power flow model in matrix form and the bus split representation. Section \ref{sec:senti} presents the linear sensitivity analysis for modeling the post-split line flows and bus phase angles. In Section \ref{sec:id}, the synchrophasor data based identification problem is formulated, along with the development of the tractable reformulation via McCormick relaxation.  Numerical studies using the IEEE 14-bus and 300-bus systems are presented in Section \ref{sec:num} to demonstrate the validity and efficiency of the proposed identification algorithm, and the paper is wrapped up in Section \ref{sec:con}.

\textit{Notation:} Upper (lower) boldface symbols stand for matrices (vectors); $(\cdot)^{\mathsf T}$ stands for transposition; $\bbI$ for identity matrix; $|\, \cdot \,|$ denotes the cardinality of a set; $\left\| \, \cdot \, \right\|$ denotes the vector norm; $\bbe_i$ denotes the standard basis vector with all entries being 0 except for the $i$-th entry equals to 1.

\section{System Modeling}\label{sec:ps}
We introduce the dc power flow model \cite{stott2009dc} for the modeling of bus split events.
{\yuqi As an approximation of nonlinear ac power flow model, the  dc model is amendable for  fast analysis with closed-form solutions. It is also known for high accuracy in approximating the real power flow, and our numerical results in Section \ref{sec:num} have corroborated the suitability of dc power flow for bus split modeling.} {\yuqi We can also extend the dc-based analysis in this work to more accurate approximation methods such as fast decoupled power flow \cite{stott1974fast} or operating-point based approximation  \cite{baldick2003variation}.}

Consider a transmission system with $(N+1)$ buses collected in the set $\cal N :=$ $\{ 0, 1,\ldots,N \}$, and $L$ transmission lines represented by the set $\cal L :=$ $\{(i,j)\}$. For each bus $i$, we use ${\cal N}_{i}$ to denote the set collecting its adjacent buses, and let $\theta_i$ denote its voltage angle, as well as $g_i$ and $d_i$ as its connected generation and load, respectively. Hence, the power injection per bus $i$ is $p_i = g_i - d_i$. Without loss of generality (Wlog), we set bus 0 to be the reference angle bus with $\theta_{0} = 0$. All non-reference phase angles are concatenated into vector \(\bm{\theta} \in \mathbb{R}^{N}\); similarly for $\mathbf{p} \in \mathbb{R}^{N}$. 
For each line $(i,j)$, the power flow from bus $i$ to bus $j$ is denoted by $f_{ij}$ and  
given by:
\begin{align}
f_{ij} =  \frac{1}{x_{ij}}(\theta_i - \theta_j) = {b_{ij}}(\theta_i - \theta_j),  \forall (i,j)\in \cal L  \label{PF}
\end{align}
where $x_{ij}$ is the line reactance and its inverse equals to $b_{ij} = {1}/{x_{ij}}$. Concatenating \eqref{PF} into matrix form gives rise to the flow vector $\mathbf{f} \in \mathbb{R}^{L}$, as:
\begin{align}
\mathbf{f} = \mathbf{K}\bm{\theta}
\label{eq:line_flow}
\end{align}
where matrix $\mathbf{K} \in \mathbb R^{L\times N}$ captures the network topology.
%
The $\ell$-th row of $\bbK$ corresponding to line $(i,j)$ is given by $
b_{ij}\left(\bbe_i-\bbe_j\right)^{\mathsf T}$, where $\bbe_i\in\mathbb R^N$ is the standard basis vector. 
Due to nodal flow conservation, one can sum up all the line flows to form the injected power $\bbp$, as given by: 
\begin{align}
\bbp = \bbB \bm{\theta}
\label{eq:DCPF2}
\end{align}
where the \textit{Bbus matrix} $\bbB\in \mathbb R^{N\times N}$ is invertible with each entry:
\begin{align}
B_{ij}=\left\{
\begin{array}{c l}	
     \sum_{k\in\ccalN_i} {b_{ik}}, & \textrm{if}\; i=j\\
     -b_{ij}, & \textrm{if}\; (i,j)\in \cal L\\
     0, & \textrm{otherwise}
\end{array}\right.
\label{eq:invertible}
\end{align}
Therefore, matrix $\bbB$ can be also given by
\begin{align}
    \bbB = \sum_{(i,j)\in\ccalL} b_{ij} (\bbe_i - \bbe_j)(\bbe_i-\bbe_j)^{\mathsf T}. \label{eq:matB}
\end{align}
By solving for $\bm{\theta}$ in \eqref{eq:DCPF2}, one can write the line flow as $\bbf = (\bbK\bbB^{-1})\bbp$, with the coefficients in $(\bbK\bbB^{-1})$ termed as the injection shift factors (ISFs) that can transform from the injection $\bbp$ to line flow $\bbf$.
%

\begin{figure}[t!]
\centering
\vspace{-2pt}
\includegraphics[trim=5.5cm 6cm 5.5cm 4.5cm,clip=true,width=.8\linewidth]{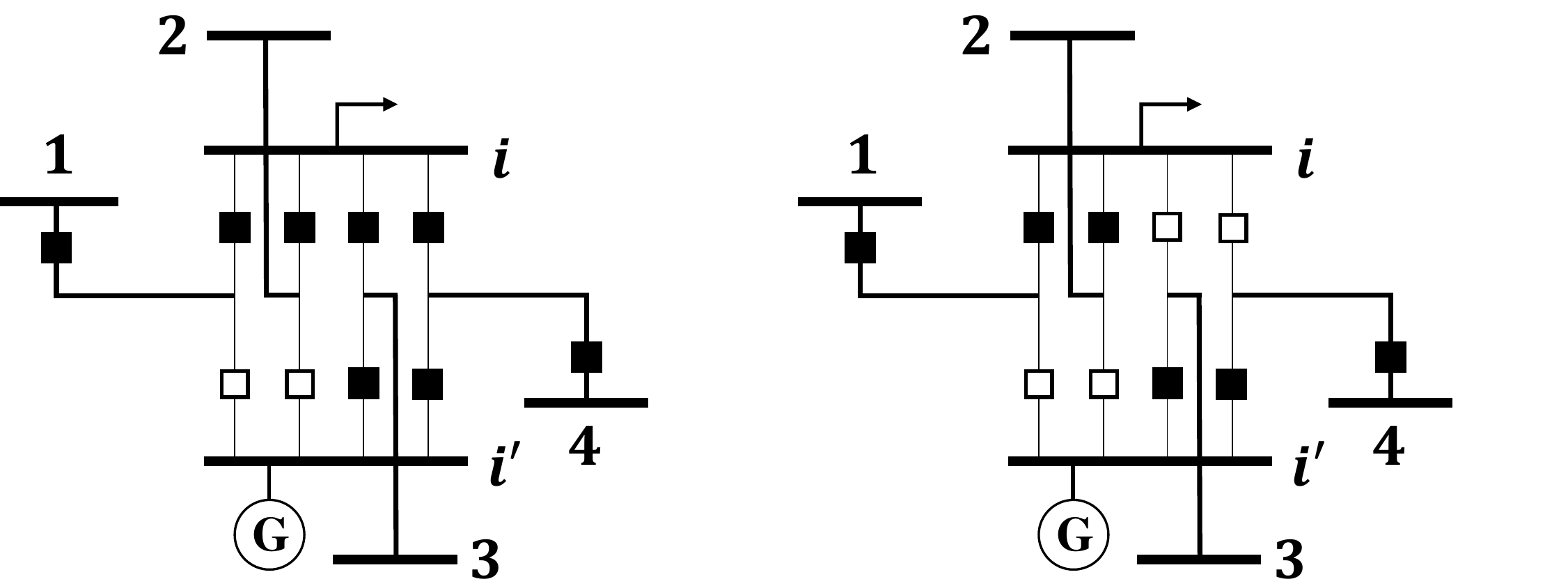}
\caption{(Left) The original substation topology and (right) the new topology with two more breakers open.}
\label{fig:bus split 2}
\end{figure}

\begin{figure}[t!]
\centering
\vspace{-2pt}
\includegraphics[trim=5.5cm 6.5cm 5cm 5cm,clip=true,width=.8\linewidth]{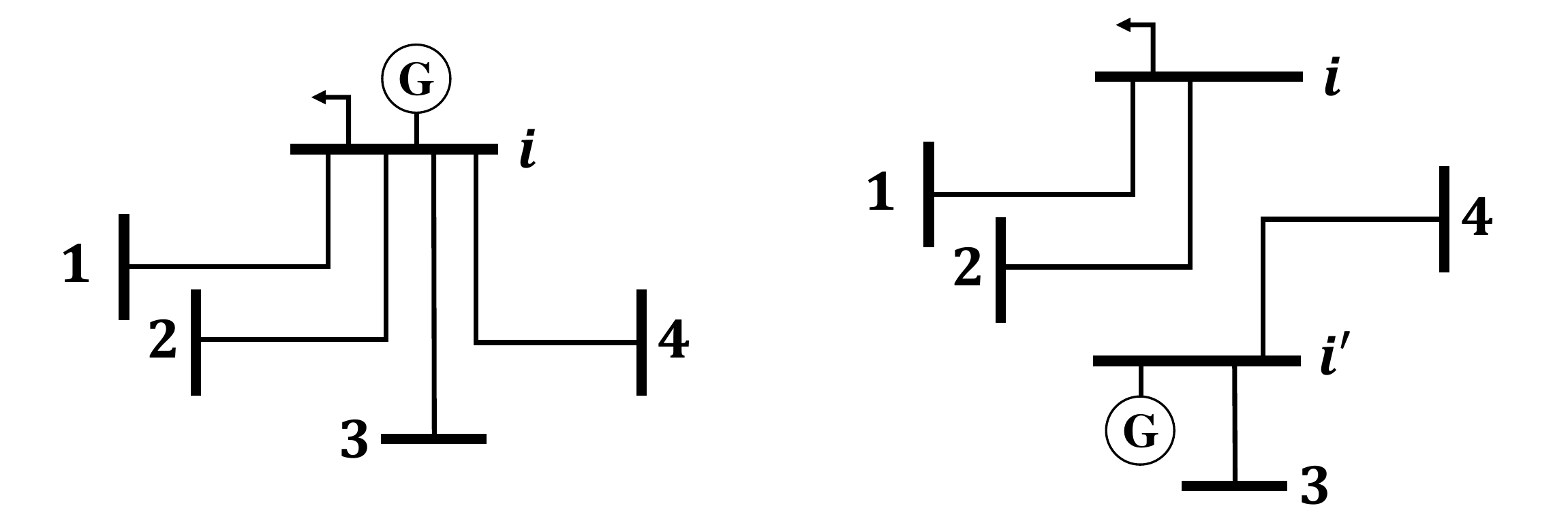}
\caption{(Left) The original bus-branch model and (right) the updated model with reconnected lines/generation/loads due to bus split.}
\vspace{-10pt}
\label{fig:bus split 1}
\end{figure}

We are interested in the grid topology changes due to \textit{bus split} within substations. Substations can be described as electrically connected nodes where multiple transmission lines terminate. 
{\yuqi 
To allow for flexible connectivity, electrical substations are equipped with switching equipment such as CBs and isolators. The latter can  disconnect or isolate different components for scheduled system maintenance or protection against faults; see e.g., \cite[Ch. 12]{santoso2018standard}. }
%
Due to CB actions, the original bus may become electrically disconnected, commonly termed as ``bus split.'' This type of topological changes is increasingly popular due to CB misoperations \cite{kekatos2012joint} or malicious cyber intrusions \cite{deka2015one,zhou2018false}. Fig. \ref{fig:bus split 2} illustrates one such event, with solid squares representing closed breakers and hollow ones for open breakers. If the top right CBs open, bus $i$ is split into two disconnected buses, $i$ and $i'$, in the same substation. Accordingly, the connectivity for transmission lines has changed, as well as for generation and load within this substation. Although the two buses are physically co-located in the same substation, they become electrically disconnected which affects the bus-branch model as shown in Fig. \ref{fig:bus split 1}. {\yuqi Note that line outages within one substation are a special case of bus split events. If the new bus $i'$ is not connected to any generation or load, the reconnected lines are equivalently in outage. Hence, bus splitting represents a general class of grid topology changes, and the ensuing sensitivity analysis and identification solution can generalize line outages as well.} 


\section{Linear Sensitivity Analysis for Bus Split}
\label{sec:senti}

\begin{figure}[t!]
\centering
\vspace{-2pt}
\includegraphics[trim=3cm 4cm 3cm 5cm,clip=true,width=.9\linewidth]{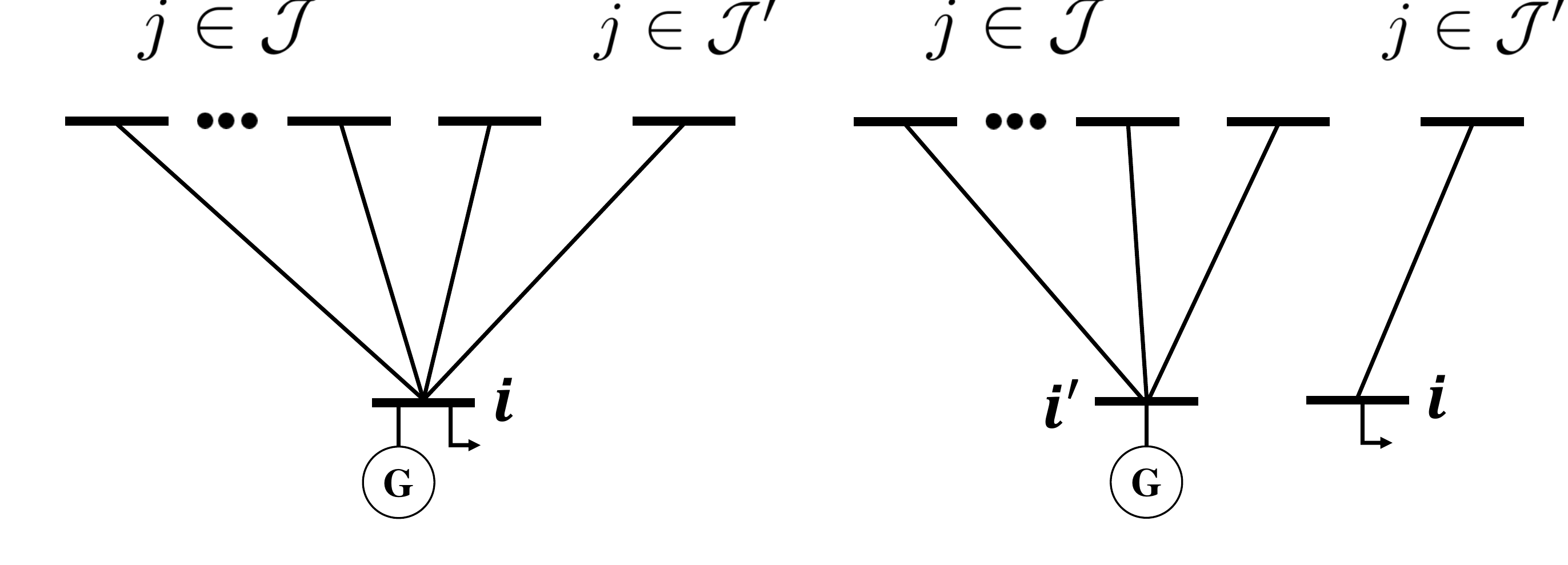}
\caption{(Left) The original bus-branch model and (right) the updated model due to bus split.} \label{fig:bus_split}
\end{figure}

Consider the split of bus $i$ leading to a new bus $i'$, as illustrated by Fig. \ref{fig:bus_split}. After the split, let the subset $\mathcal{J} \subseteq \mathcal{N}_{i}$ collect the adjacent buses reconnected to the new bus $i'$. Moreover, denote $\tdbbtheta = \left[\bbtheta';~\tdtheta_{i'}\right]\in \mathbb{R}^{N+1}$ as the new angle vector after splitting; and similarly for $\tdbbp$ and ${\tilde{\mathbf{B}}}$. The vector/matrix dimension increases by one due to the new bus $i'$. The goal is to update the new angle $\tdbbtheta$ and line flow $\tdbbf$. 
{\yuqi
Bus split events may disconnect the whole network and cause islanding in an interconnected system \cite{sun2003splitting}. Nonetheless, islanding arises very infrequently due to the meshed topology of transmission systems \cite{kyriacou2017controlled}. Even if system islanding happens, the lack of synchronization among islands leads to very noticeable response signatures such as separated frequency or phase angle. For simplicity, this paper  assumes no island is formed from the split of bus $i$.}

Matrix ${\tilde{\mathbf{B}}}$ follows a block structure as:
\begin{align}
{\tilde{\mathbf{B}}} := 
\begin{bmatrix} 
\mathbf{B'} & \bm{\ell} \\
{\bm{\ell}}^{\mathsf T} & d 
\end{bmatrix}
\label{eq:B_tilde}
\end{align}
where $\mathbf{B}'$ is the $N  \times N$ sub-matrix, while vector $\bm{\ell}$ and scalar $d$ capture the rest. Based on Fig. \ref{fig:bus_split}, the lines $(i,j)$ with $j \in \mathcal{J}$ denote those lines that are reconnected to the new bus $i'$ after the split. Similarly, we use buses $j \in \mathcal{J}'$ to match those lines $(i,j)$ that remain connected to bus $i$ after the split. Hence, the submatrix $\mathbf{B}'$ is formed by eliminating these reconnected lines from $\bbB$, as given by [cf. \eqref{eq:matB}]
\begin{align}
\bbB' = \bbB  - \sum_{j \in \mathcal{J}} b_{i{j}}\left[{(\mathbf{e}_{i} - \mathbf{e}_{j})}{(\mathbf{e}_{i} - \mathbf{e}_{j})}^{\mathsf T} - \mathbf{e}_{j} \mathbf{e}_{j}^{\mathsf T}\right].
\label{eq:b_prime}
\end{align}
The reconnected lines $\{(i',j)\}$ also affect the rest of $\tdbbB$, as 
\begin{align}
\bm{\ell} =  -\sum_{j \in \mathcal{J}} b_{i j}\mathbf{e}_{j}, \quad \textrm{and} \quad  d = \sum_{j \in \mathcal{J}} b_{i j}.
\label{eq:ld}
\end{align}
With these definitions, one can find the inverse of ${\tilde{\mathbf{B}}}$ using the popular \textit{matrix inverse lemma} \cite[p. 650]{boyd2004convex}. To this end, define the following \textit{Schur complement} of the entry $d$ to matrix $\tdbbB$ as  
\begin{align}
\bbB_d &  \coloneqq \mathbf{B'} -   \left(d^{-1}\right) \bm{\ell}{\bm{\ell}}^{\mathsf T}. 
\label{eq:Sd}
\end{align}
Using $\bbB_d$, we can obtain the block structure of the inverse as
\begin{align}
{\tilde{\mathbf{B}}}^{-1} = \begin{bmatrix} 
\bbB_d^{-1} & - \left(d^{-1}\right)\bbB_d^{-1}\bm{\ell} \\
-\left(d^{-1}\right){\bm{\ell}}^{\mathsf T}\bbB_d^{-1} & 
({d}^{-2}){\bm{\ell}}^{\mathsf T}\bbB_d^{-1}{\bm{\ell}}+ {d}^{-1}
\end{bmatrix}.
\label{eq:blockwise3}
\end{align}
Note that one can rewrite $\bbB_d = \bbB - d(\bbu \bbu^{\mathsf T})$ with the vector
\begin{align}
\bbu : = \mathbf{e}_{i} + \left(d^{-1}\right)\bm{\ell}.
\end{align}
For example, if $|\mathcal{J}| = 1$, it becomes the case of  single line reconnection as discussed in  \cite{zhou2019bus} where vector $\bbu = \mathbf{e}_{i} - \mathbf{e}_{j}$ with $\ccalJ = \{j\}$. Applying the Sherman-Morrison formula \cite{sherman1950adjustment} leads to the inverse
\begin{align}
\bbB_d^{-1} 
& = (\bbI +\bbDelta_\bbu) \bbB^{-1} = \mathbf{B}^{-1} + \frac{d \mathbf{B}^{-1}\bbu\bbu^{\mathsf T}}{1-d\bbu^{\mathsf T}\mathbf{B}^{-1}\bbu}\mathbf{B}^{-1}, \label{eq:Sdinv}
\end{align}
with $\bbDelta_\bbu$ capturing the fractional term above. Clearly, the inverse $\bbB_d^{-1} $ can be quickly formed with $\bbB^{-1}$ available, and so are the other blocks in \eqref{eq:blockwise3}.

The next step is to form the post-split power injection vector $\tdbbp$. Fig. \ref{fig:bus_split} indicates that generation and load originally connected to bus $i$ could become attached to the new bus $i'$ due to the split. Let $\tdp_i$ denote the total injected power to bus $i'$ after the split, and thus the injected power to bus $i$ becomes $(p_i - \tdp_i)$. For example, Fig. \ref{fig:bus_split} shows the case of $\tdp_i = g_i$ and  $(p_i - \tdp_i) = -d_i$. In general, $\tdp_i$ can also be the attached load, or, the combination of generation and load. Therefore, we can express the injected power as 
\begin{align}
\tilde{\mathbf{p}} = 
\begin{bmatrix} 
\mathbf{p}\\
0
\end{bmatrix} + 
\begin{bmatrix} 
- {\tdp_i}\mathbf{e}_i\\
{\tdp_i}
\end{bmatrix} = 
\begin{bmatrix} 
\mathbf{p} - {\tdp_i}\mathbf{e}_i\\
{\tdp_i}
\end{bmatrix}.
\label{eq:p_split}
\end{align}

Using \eqref{eq:blockwise3} and \eqref{eq:p_split}, one can obtain the post-split  $\tilde{\bm{\theta}}$ as 
\begin{align}
\begin{bmatrix} 
\bm{\theta}'\\ \tdtheta_{i'}\end{bmatrix} =  {\tilde{\mathbf{B}}}^{-1} \tilde{\mathbf{p}} =
\begin{bmatrix} 
\bm{\theta}'\\
 d^{-1} \left(-\bm{\ell}^{\mathsf T} \bm{\theta}'+
{\tdp_i}\right)
\end{bmatrix} 
\label{eq:theta_tilde}
\end{align}
where the angle vector for the original $N$ buses is
\begin{align}
    \bm{\theta}'\coloneqq \bbB_d^{-1} \left(\mathbf{p} - \tdp_i\mathbf{e}_i - d^{-1} \tdp_i \bm{\ell}  \right)
    \label{eq:theta_prime}
\end{align}

\begin{figure}[t!]
\centering
\vspace{-2pt}
\includegraphics[trim=3cm 4cm 3cm 4cm,clip=true,width=.9\linewidth]{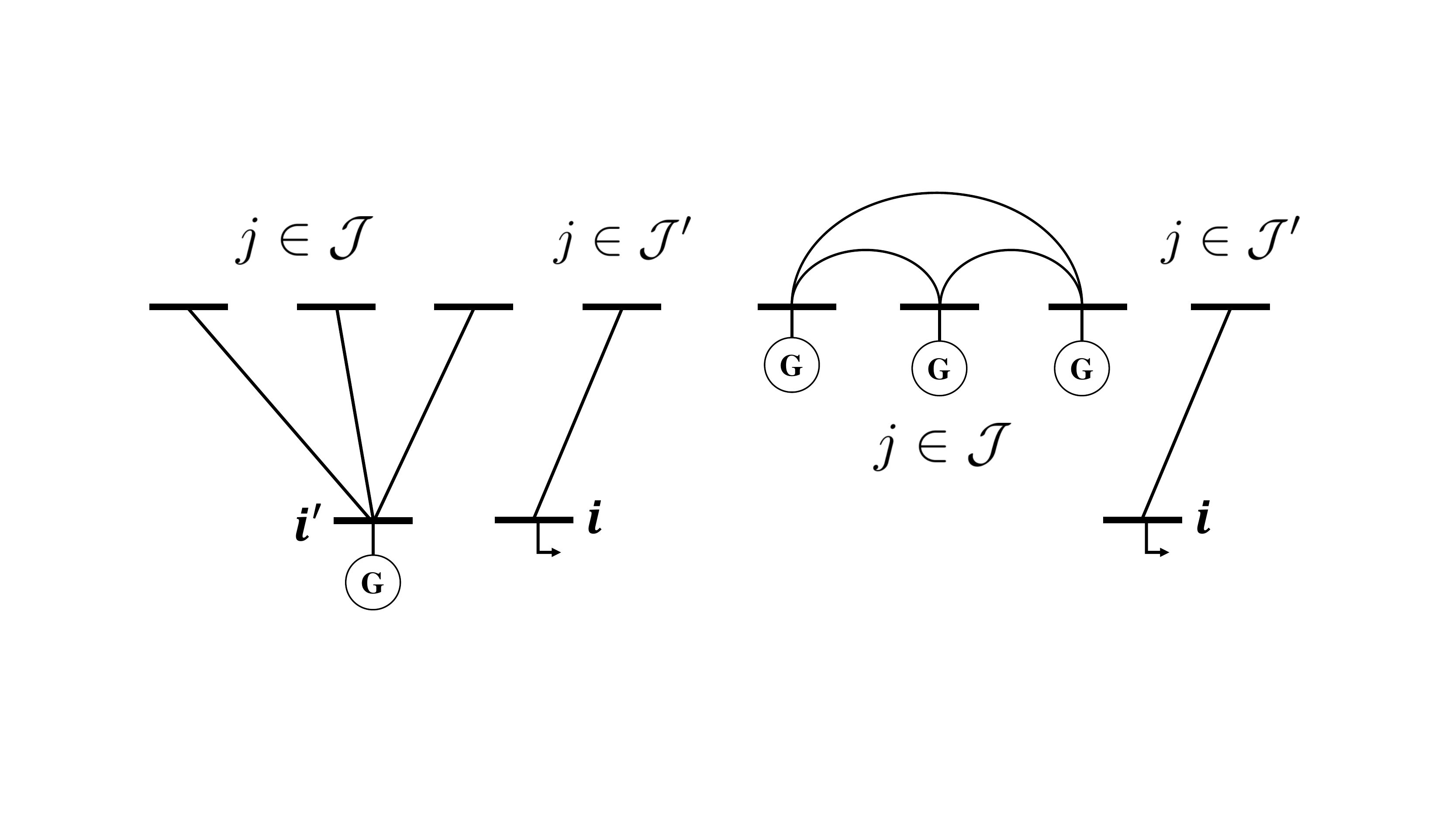}
\caption{(Right) The equivalent bus-branch model after eliminating bus $i'$ for (left) the system due to the bus split.} \label{fig:reduction}
\end{figure}

Although the solution $\tdbbtheta$ can be readily computed once obtaining  \eqref{eq:Sdinv}, it is possible to better interpret it by developing an equivalent model for the post-split system. As illustrated in Fig. \ref{fig:reduction}, the new bus $i'$ can be eliminated by connecting every pair of buses in $\mathcal{J}$. For example, the new line connecting any buses $j_1, j_2\in \mathcal{J}$ would have the equivalent line parameter $b_{j_1 j_2}:= \left(d^{-1}\right)b_{i j_1} b_{i j_2}$ by recalling $d$ from \eqref{eq:ld}. This follows from the well-known Kron's reduction approach \cite{dorfler2012kron}. Similarly, after eliminating bus $i'$, its injected power $\tdp_i$ is reallocated to each of its neighboring buses, with $d^{-1} {b_{i j}} \tdp_i$ to bus $j \in \mathcal{J}$.

Interestingly, the equivalent topology in Fig. \ref{fig:reduction}  exactly matches the matrix $\bbB_d$ as defined in \eqref{eq:Sd}. Substituting \eqref{eq:b_prime} into \eqref{eq:Sd}, one can decompose $\bbB_d = \mathbf{B}^{\ccalJ} + \mathbf{B}^{a}$, where the matrix
\begin{align}
\mathbf{B}^{\ccalJ} := \bbB  - \sum_{j \in \mathcal{J}} b_{i{j}}{(\mathbf{e}_{i} - \mathbf{e}_{j})}{(\mathbf{e}_{i} - \mathbf{e}_{j})}^{\mathsf T} \label{eq:BJ}    
\end{align} 
corresponds to the Bbus matrix with the outage of lines in $\{(i,j)\}_{j \in \mathcal{J}}$ from the original system, while the remaining part  
\begin{align}
{\mathbf{B}}^{a} = \sum_{j \in \mathcal{J}} b_{i j} \mathbf{e}_{j} \mathbf{e}_{j}^{\mathsf T} -\left(d^{-1}\right) \bm{\ell}{\bm{\ell}}^{\mathsf T}
\label{eq:B_a}
\end{align}
reflects the new equivalent lines among all buses in $\mathcal{J}$. 
Hence, the angle vector $\bbtheta'$ in \eqref{eq:theta_prime} perfectly matches the equivalent system for the original $N$ buses with updated topology/injection, as its total injection now includes the new ones at each $j\in\ccalJ$ due to eliminating bus $i'$. 
Using the updated angles in $\boldsymbol{\theta}'$, one can further recover the angle at bus $i'$ by solving the total power flow on the eliminated lines in $\{(i',j), ~\forall j \in \mathcal{J}\}$ as
\begin{align}
{\tdp_i}= \sum_{j \in \mathcal{J}} b_{i j} \left(\tdtheta_{i'} - \mathbf{e}_j^{\mathsf T}\bm{\theta}'\right),
\label{eq:last_entry}
\end{align}
leading to the same solution of $\tdtheta_{i'}$ as in \eqref{eq:theta_tilde}. This equivalencing analysis clearly explains how the post-split system can be related to the original one. {\yuqi Interestingly, this analysis also includes the case of multiple line outages within substation $i$, since the latter is a special case of bus splitting.}

Upon solving for the post-split angle $\tilde{\bbtheta}$, one can formalize the sensitivity analysis in terms of the changes of bus angles and line flows, as given in the following proposition.
\begin{proposition}
\label{prop:1}
For the split of bus $i$ with a re-connected injection of $\tdp_i$, the phase angle difference at the original $N$ buses can be written as
\begin{align}
\bbtheta'-\bm{\theta} =
\bbDelta_\bbu \bbtheta - (\bbI+\bbDelta_\bbu)\bbB^{-1}\left( \tdp_i\mathbf{e}_i + d^{-1} \tdp_i \bm{\ell}\right)
\label{eq:delta_theta}
\end{align}
where $d$, $\bm{\ell}$, and  $\bbDelta_\bbu$ are given in \eqref{eq:ld} and \eqref{eq:Sdinv}.
Accordingly, the line flow difference due to the bus split is [cf. \eqref{eq:line_flow}]
\begin{align}
\tdbbf - \mathbf{f} =
\mathbf{K}\bbDelta_\bbu  \bbtheta - \mathbf{K}(\bbI+\bbDelta_\bbu)\bbB^{-1}\left( \tdp_i\mathbf{e}_i + d^{-1} \tdp_i \bm{\ell}\right).
\label{eq:delta_f}
\end{align}
\end{proposition}
Strictly speaking, \eqref{eq:delta_f} holds for the lines other than the reconnected ones $\{(i',j)\}_\{j \in \mathcal{J}\}$, as the flows of the latter lines are related to $\tdtheta_{i'}$.
In addition, since matrix $\bbDelta_\bbu$ is a linear transformation of $\bbB^{-1}$ [cf. \eqref{eq:Sdinv}],  the matrix product $\mathbf{K}\bbDelta_\bbu$ in \eqref{eq:delta_f} can be quickly computed using the ISF matrix $\bbK\bbB^{-1}$.

Given the detailed post-split topology and injection, the sensitivity analysis in Proposition \ref{prop:1} enables explicit expression of the phase angle and line flow change. In addition to generalizing line outage sensitivity analysis, this result could benefit security-constrained economic dispatch by including bus split contingencies   \cite{zhou2019bus} and other related EMS tasks. 
The sensitivity analysis results could also be useful for the monitoring of bus split events, which is the subject of the ensuing section. 

\section{Tractable Identification via McCormick Relaxation}
\label{sec:id}

Increasing deployment of high-resolution sensors such as PMUs allows for real-time identification of anomalous events such as topology changes. It is thus possible to utilize synchrophasor data to efficiently identify the underlying CB statuses during bus split events. Based on the sensitivity analysis in Sec. \ref{sec:senti}, an intuitive solution could be the exhaustive search by enumerating all number of possible scenarios. However, the complexity order of this approach is not ideal for real-time identification and motivates us to consider more efficient alternatives.  


As in the identification of line outages \cite{tate2008line,zhu2012sparse}, synchrophasor data provides the difference between pre- and post-event phase angle measurements which can be used to locate bus splits. 
Consider again the split of bus $i$ as in Sec. \ref{sec:senti}. 
This candidate bus can be selected based on the bus locations where significant angle changes have been observed, as detailed in Sec. \ref{sec:num}.
Given bus $i$ and its subset $\ccalJ$ that include all reconnected buses/generation/loads, the post-split phase angle $\tilde{\bm{\theta}}$ can be computed in \eqref{eq:theta_tilde}. This suggests a brute-force solution by enumerating all possible scenarios of sets $\ccalJ$ and selecting the one best matching the measured changes. Towards this end, for the target bus $i$ let us define the binary variables $\{z^j\}_{j \in \mathcal{N}_{i}}$ to indicate the post-event status for each line $(i,j)$, such that 
\begin{align}
z^j = \left\{
\begin{array}{c l}	
     1, & \textrm{if bus $j$ re-connects to bus $i'$}\\
     0, & \textrm{if bus $j$ stays connected to bus $i$}
\end{array}\right.
\end{align}
Similarly, let the binary variables $\{z^g,~z^d\}$ denote the post-event connection status for the generation and load located at bus $i$. For simplicity, only one generator and one load are considered here, which can be extended to the case of multiple number of them. Given all these binary variables in vector $\bm z\in \{0,1\}^{({|\ccalN_i|} + 2)}$, one can use the sensitivity analysis to compute the post-split angle $\tilde{\bm{\theta}}(\bm z)$. Denoting the measured post-split angle as vector $\bbc$, the identification problem becomes to minimize the norm of mismatch error between the two, namely $\| \tilde{\bm{\theta}}(\bm z) - \bbc \|$. To address the potential approximation error of dc power flow model, a more popular error objective is to compare the \textit{phase angle change} instead of the post-split phase angle itself; see e.g., \cite{zhu2012sparse}. Hence, the bus split identification problem becomes
\begin{align}
\min_{\bm z \in \{0,1\}^{({|\ccalN_i|} + 2)}} \quad & \left\| \left( \tilde{\bm{\theta}}(\bm z) - \hat{\bm{\theta}} \right)  - \boldsymbol \delta \right\|
\label{eq:enumeration_1}
\end{align}
where $\hat{\bm{\theta}} \coloneqq \left[\bbtheta;\theta_{i} \right]$ is pre-split phase angle solution by the dc model and $\boldsymbol \delta$ is the observed phase angle difference by the PMUs. Common choice of norms such as $L_1$, $L_2$ or $L_{\infty}$ can be used for quantifying the mismatch error in the objective \eqref{eq:enumeration_1}. In this paper, we select the $L_1$ norm error which is known to enjoy a nice reformulation to linear objective function.
Note that if only partial angles are monitored or the line flow data is also available, one can modify the objective function of \eqref{eq:enumeration_1} by extracting corresponding entries or adding a linear transformation of the angle difference vector. Such generalization will be investigated numerically in Sec. \ref{sec:num_300}. Using \eqref{eq:enumeration_1}, the naive exhaustive search approach works by enumerating all $2^{({|\ccalN_i|} + 2)}$ scenarios of $\bm z$, and picking the one with the minimum error in \eqref{eq:enumeration_1}. Such exponential complexity is not suitable for real-time implementation, and will be addressed here by developing a tractable solution. {\yuqi Although simultaneous splits at multiple locations are much less likely to occur in view of the fast rates of synchrophasor samples, identification of such events may be included by \eqref{eq:enumeration_1} as well through introducing binary variables for all candidate substations.}



The main challenge in directly optimizing for $\bm{z}$ in \eqref{eq:enumeration_1} lies in the post-split power flow model, as
\begin{align}
\tilde{\bm{\theta}}(\bm z) = \left[\tilde{\mathbf{B}}(\bm{z})\right]^{-1} \tilde{\mathbf{p}}(\bm{z}).
\label{eq:pf_z}
\end{align}
The post-split Bbus matrix depends on the line status variables $\{z^j\}$, given by
\begin{align}
\tilde{\mathbf{B}}(\bm{z}) & = \bar{\mathbf{B}} + \sum_{j \in \mathcal{N}_{i}} z^{j} ~ \mathbf{D}^{j} \label{eq:Bz}
\end{align}
where $\bar{\mathbf{B}} = \left[\mathbf{B} \; \bm{0}; \bm{0}^{\mathsf T}\right]$ corresponds to the original Bbus matrix for the augmented system including the new bus $i'$, while the change due to each reconnected line is captured by
\begin{align}
\mathbf{D}^{j}\!:= b_{ij} \left[{(\mathbf{e}_{i'} - \mathbf{e}_{j})}{(\mathbf{e}_{i'} - \mathbf{e}_{j})}^{\mathsf T} \! -\! (\mathbf{e}_{i} - \mathbf{e}_{j})(\mathbf{e}_{i} - \mathbf{e}_{j})^{\mathsf T}\right]\!.\label{eq:Dj}
\end{align}
Similarly, the post-split injection vector depends on the generation/load status variables, as 
\begin{align}
\tilde{\mathbf{p}}(\bm{z}) = \bar{\mathbf{p}} + z^g\boldsymbol{\delta}^{{g}} - z^d\boldsymbol{\delta}^{{d}} 
\label{eq:p_tilde}
\end{align}
where $\bar{\mathbf{p}} = \left[\mathbf{p}; 0\right]$ also augments the dimension of $\mathbf{p}$ to include the new bus $i'$, while the other two vectors are given by
\begin{align}
    \label{eq:deltagd}
\boldsymbol{\delta}^{{g}} := g_i (\mathbf{e}_{i'} - \mathbf{e}_{i}) \textrm{~and~} \boldsymbol{\delta}^{{d}}:= d_i (\mathbf{e}_{i'} - \mathbf{e}_{i}).
\end{align}

Substituting \eqref{eq:pf_z} into the  objective of \eqref{eq:enumeration_1} introduces the term $[\tilde{\mathbf{B}}(\bm{z})]^{-1}$, which is the inverse of a matrix function of binary variables.  Furthermore, its multiplication with $\tilde{\mathbf{p}}(\bm{z})$ poses additional nonlinearity (bilinearity) to the problem. Both issues lead to the lack of tractability in solving \eqref{eq:enumeration_1}.

To tackle these issues, we propose to adopt the \emph{McCormick relaxation} technique \cite{McCormick1976computability}, which is a powerful tool for dealing with bilinear terms \cite{gupte2013solving}. {\yuqi It converts the latter to linear constraints that can are provably \textit{equivalent} for our problem.} Specifically, by introducing an additional matrix $\mathbf{X} \in \mathbb R^{(N+1) \times (N+1)}$ to represent $[\tilde{\mathbf{B}}(\bm{z})]^{-1}$, we can reformulate \eqref{eq:enumeration_1} as: 
\begin{subequations} \label{eq:mc_1}
\begin{align}
\min_{\bm{z}, \mathbf{X}} \quad & \left\|  {\mathbf{X}} \tilde{\mathbf{p}}(\bm{z}) - \hat{\bm{\theta}} - \boldsymbol \delta \right\| \label{eq:mc_1_a}\\
\textrm{s.t.} \quad & \tilde{\mathbf{B}}(\bm{z}) \mathbf{X} = \mathbf{I}. \label{eq:mc_1_b}
\end{align}
\end{subequations}
This way, the matrix inversion is no longer needed, which is replaced by the bilinear products between the unknowns.
Note that the bilinear constraint \eqref{eq:mc_1_b} for enforcing the relation between $\bbX$ and $\tdbbB(\bm{z})$ now becomes [cf. \ref{eq:Bz}]:
\begin{align}
\bar{\mathbf{B}}\mathbf{X} + \sum_{j \in \mathcal{N}_{i}}  \mathbf{D}^{j} \mathbf{Y}^j = \mathbf{I}
\label{eq:eq_1}
\end{align}
with the product $\mathbf{Y}^{j} :=z^{j} \mathbf{X}$ defined for each line $(i,j)$. Given the binary $z_j\in\{0,1\}$ and the $(m,n)$-th entry $X_{mn} \in [{X}^{\min}_{mn},{X}^{\max}_{mn}]$, each entry of $\mathbf{Y}^{j}$ is written as:
\begin{align}
{Y}^{j}_{mn} = z^{j} {X}_{mn},~\forall~m,n.
\label{eq:equality_1}
\end{align}
{\yuqi We can show that the bilinear relation in \eqref{eq:equality_1} is equivalent to the following  four \textit{linear} inequalities:}
\begin{subequations} \label{eq:mc}
\begin{align}
& {Y}^{j}_{mn} \geq z^{j} {X}^{\min}_{mn},\label{eq:mc_c}\\
& {Y}^{j}_{mn} \geq {X}_{mn} + z^{j} {X}^{\max}_{mn} - {X}^{\max}_{mn},\label{eq:mc_d}\\
& {Y}^{j}_{mn} \leq z^{j} {X}^{\max}_{mn},\label{eq:mc_e}\\
& {Y}^{j}_{mn} \leq {X}_{mn} + z^{j} {X}^{\min}_{mn} - {X}^{\min}_{mn}. \label{eq:mc_f}
\end{align}
\end{subequations}
{\yuqi 
Each inequality in \eqref{eq:mc} can be verified by {\yuqi considering the upper/lower bounds of $z_j$ and ${X}_{mn}$}. Interestingly, the set of four inequalities in \eqref{eq:mc} jointly guarantees that \eqref{eq:equality_1} would hold for any binary $z^j$. To demonstrate this, first consider the case of $z^{j} = 0$. Constraints \eqref{eq:mc_c} and \eqref{eq:mc_e} together lead to  $Y_{mn}^{j} = 0$, and thus \eqref{eq:equality_1} holds. Otherwise if $z^{j} = 1$, the other constraints \eqref{eq:mc_d} and \eqref{eq:mc_f} would enforce that $Y_{mn}^{j} = {X}_{mn}$. Hence, for binary $z^j$ the set of inequalities in \eqref{eq:mc} is \textit{equivalent to} the bilinear relation in \eqref{eq:equality_1}. 
}
Reformulating \eqref{eq:equality_1} using the linear inequalities in \eqref{eq:mc} is known as the \textit{Mccormick relaxation} technique and has been popularly used in power system topology  analysis and design  \cite{kocuk2017new,bhela2019designing,park2020optimal}. 

For the objective function in \eqref{eq:mc_1_a}, we can introduce matrices $\bbW^p := z^p \bbX$ with $p \in \{g,d\}$ indicating either the generation or load, and similarly convert the matrix products into equivalent linear inequalities. We replace the resultant bilinear terms in \eqref{eq:mc_1_a} for $p \in \{g,d\}$ using the follows:
\begin{subequations}\label{eq:mcp}
\begin{align}
& {W}^{p}_{mn} \geq z^{p} {X}^{\min}_{mn}, \label{eq:mc_3_aa}\\
& {W}^{p}_{mn} \geq {X}_{mn} + z^{p} {X}^{\max}_{mn} - {X}^{\max}_{mn}, \label{eq:mc_3_bb}\\
& {W}^{p}_{mn} \leq z^{p} {X}^{\max}_{mn}, \label{eq:mc_3_cc}\\
& {W}^{p}_{mn} \leq {X}_{mn} + z^{p} {X}^{\min}_{mn} - {X}^{\min}_{mn}. \label{eq:mc_3_dd}
\end{align}
\end{subequations}
\begin{proposition}
\label{prop:2}
Given the bounds $[{X}^{\min}_{mn},{X}^{\max}_{mn}]$ for each entry ${X}_{mn}$, \eqref{eq:mc_1} {\yuqi is equivalent to} the following mixed-integer programming problem upon substituting $\tilde{\mathbf{p}}(\bm{z})$ into \eqref{eq:p_tilde}:
\begin{subequations} \label{eq:mc_3}
\begin{align}
\min_{\bm{z}, \mathbf{X}, \{\mathbf{Y}^j\}, \{\mathbf{W}^p\}} \quad & \| {\mathbf{X}}\bar{\mathbf{p}} + \mathbf{W}^g \boldsymbol{\delta}^{{g}} - \mathbf{W}^d \boldsymbol{\delta}^{{d}} - \hat{\bm{\theta}} - \boldsymbol \delta \| \label{eq:mc_3_a}\\
\textrm{s.t.} \quad & \eqref{eq:eq_1},~
\eqref{eq:mc}\mathrm{~and~}\eqref{eq:mcp},~\forall~(m,n).
\label{eq:2mc}
\end{align}
\end{subequations}
\end{proposition}
As mentioned earlier, the $L_1$ error norm in \eqref{eq:mc_3_a} can lead to an equivalent linear objective cost; see e.g., \cite[Ch. 4]{boyd2004convex}. Meanwhile, all constraints in \eqref{eq:2mc} are linear, thanks to the Mccormick relaxation technique.
Hence, the original problem  \eqref{eq:enumeration_1} is converted to a mixed-integer linear program (MILP), for which there exist several off-the-shelf efficient solvers such as CPLEX and Gurobi. 
This equivalent MILP reformulation \eqref{eq:mc_3} constitutes as a tractable solution for identifying the status of multiple connections within a specific substation.  The setting of tight bounds $[{X}^{\min}_{mn},{X}^{\max}_{mn}]$ for each system will be discussed soon. 

It is possible to further reduce the number of decision variables in \eqref{eq:mc_3}  by taking advantage of the sparse structure of matrix $\mathbf{D}^{j}$. To this end, let us define the product $\boldsymbol{{\Delta}}^{j}:= \mathbf{D}^{j} \mathbf{Y}^j = z^j \mathbf{D}^j \mathbf{X}$, which was used in \eqref{eq:eq_1}. With the definition of $\mathbf{D}^j$ in  \eqref{eq:Dj}, the product $\boldsymbol{{\Delta}}^{j}$ would involve only three rows of  $\mathbf{X}$, namely the $i$-th, $j$-th and $i'$-th rows. Accordingly, it suffices to form the corresponding submatrix $\mathbf{Y}^{j} \in \mathbb R^{3 \times (N+1)}$  as decision variables with $\forall n = 1,\ldots, N+1$: 
\begin{subequations} \label{eq:sparse_1}
\begin{align}
& Y_{1,n}^{j} = z^{j} X_{i,n}~,  \\
& Y_{2,n}^{j} = z^{j} X_{j,n}~,\\
& Y_{3,n}^{j} = z^{j} X_{i',n}~.
\end{align}
\end{subequations}
As a result, matrix $\boldsymbol{{\Delta}}^{j}$ is a sparse matrix with all nonzero elements listed here:
\begin{subequations} \label{eq:sparse_2}
\begin{align}
& \boldsymbol{{\Delta}}^{j}_{i,n} = b_{ij} \left(Y^{j}_{2,n} - Y^{j}_{1,n} \right), \\
& \boldsymbol{{\Delta}}^{j}_{j,n} = b_{ij} \left(Y^{j}_{1,n} - Y^{j}_{3,n} \right), \\
& \boldsymbol{{\Delta}}^{j}_{i',n} = b_{ij} \left(Y^{j}_{3,n} - Y^{j}_{2,n} \right).
\end{align}
\end{subequations}
The same technique can be applied to reduce the number of variables in the product of $\mathbf{W}^p$ and $\boldsymbol{\delta}^{p}$ in \eqref{eq:mc_3_a}. For the example of generation, due to the definition of $\boldsymbol{\delta}^g$ in \eqref{eq:deltagd}, the operation only depends on the $i$-th and $i'$-th columns of $\bbX$. Hence, one only needs to introduce a submatrix  $\mathbf{W}^{g} \in \mathbb R^{(N+1) \times 2}$ as decision variables and form the product similar to \eqref{eq:sparse_2}. Thanks to the sparsity of matrix/vector in \eqref{eq:mc_3}, we can reduce the number of decision variables from $\mathcal O(N^2)$ to $\mathcal O(N)$, and same for  the number of linear inequality constraints in \eqref{eq:2mc}. This simplification step will allow for achieving the identification solution more efficiently by the MILP solvers.

\color{black}





\subsection{Bounds on Continuous Variables $X_{mn}$}

With the MILP reformulation given in \eqref{eq:mc_3}, the problem remains to derive the upper/lower bounds of continuous variables in $\bbX$.
In general, the tighter the bounds are, the faster the MILP can be solved; see e.g., \cite{ostrowski2011tight}.
Hence, we will aim to find reasonably good bounds for entries of $\bbX$.

To this end, first notice that each entry $X_{mn}$ is non-negative; i.e., the lower bound $X_{mn}^{\min} = 0$ holds.
This is because the Bbus matrix $\tilde{\mathbf{B}}({\bm{z}})$ is an M-matrix \cite{plemmons1977m} with its off-diagonal entries being non-positive and real eigenvalues being non-negative, and thus the entries of its inverse are non-negative. To obtain the upper bound ${X}^{\max}_{mn}$, it is possible to directly maximize the value of $X_{mn}$ under the same constraints in the optimization problem \eqref{eq:mc_3}, as given by
\begin{subequations} \label{eq:mc_5}
\begin{align}
\max_{\bm{z}, \mathbf{X},  \{\mathbf{Y}^j\}, \{\mathbf{W}^p\}, X^{\max}} \quad & {X}^{\max} \\
\textrm{s.t.} \quad &  X_{mn} \leq {X}^{\max},~\forall (m,n) \\
\quad & \eqref{eq:eq_1},~
\eqref{eq:mc}\mathrm{~and~}\eqref{eq:mcp},~\forall~(m,n).
\end{align}
\end{subequations}

Compared with \eqref{eq:mc_3}, this new problem \eqref{eq:mc_5} differs only in the objective function and it is an MILP as well.
Clearly, its optimum solution of $X^{\max}$ defines an upper bound for every entry of $\bbX$.
Notice that to solve \eqref{eq:mc_5}, 
we need to start with a rough estimate of the upper bound of $X_{mn}$ to tackle the inequality constraints in \eqref{eq:mc} and \eqref{eq:mcp}. This is possible by using a very large value of the upper bound estimate which is refined by \eqref{eq:mc_5}. Since this upper bound holds for any choice of $\bm{z}$, this refinement step can be performed off-line and does not affect the real-time identification time.

\color{black}

\section{Numerical Results}
\label{sec:num}

In this section, we use IEEE 14-bus and 300-bus test cases to demonstrate the identification of bus split events. 
The 14-bus case allows to better illustrate the system-wide effects of bus split, while the 300-bus case is used to provide quantifiable identification error performance. 
{\yuqi The optimization problems have been implemented on a regular laptop equipped with Intel\textsuperscript{\textregistered} CPU @ 2.60 GHz and 12 GB of RAM using the MATLAB\textsuperscript{\textregistered} R2018a simulator. The power system analysis has been simulated using the ac power flow solver of MATPOWER, and the reformulated MILP-based identification problems have been solved by Gurobi. The solver was set up to utilize up to 8 available threads with a solution tolerance of $1\mathrm{e}{-04}$.
}

\subsection{{IEEE} 14-Bus System Tests}
\label{sec:num_14}

The IEEE 14-bus system consists of 20 lines and 5 conventional generators. The proposed identification algorithm has been shown effective for all possible bus split scenarios (excluding the islanding ones) for this system. All types of bus injections, generation only, load only, or the combination of both, have been considered as well. 
Instead of listing all the results, we pick the split of bus $i=13$ as an illustrative example. The first neighbor buses of bus $i=13$ is given by $\mathcal{N}_{i} = \{6, 12, 14\}$. The line $(13,14)$ and load $d_i$ are reconnected to the new bus $i'=15$ after the bus split. 

\begin{figure}[t!]
\centering
\vspace{-2pt}
\includegraphics[trim=0cm 0cm 1cm 0cm,clip=true,totalheight=0.16\textheight]{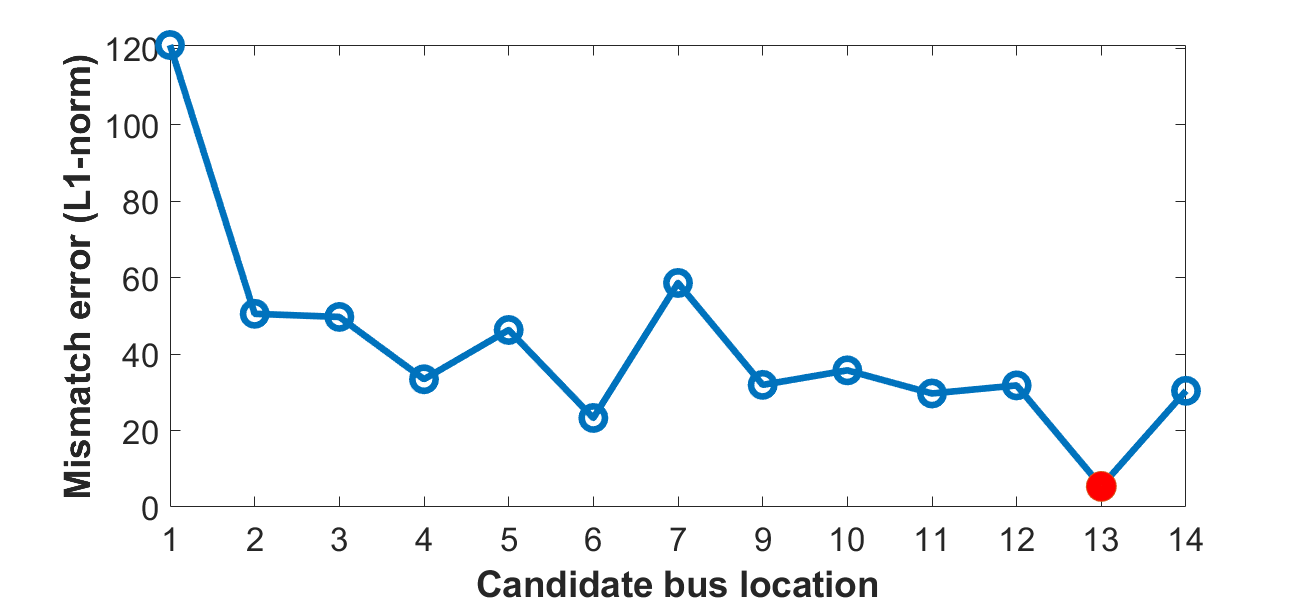}
\caption{Comparison of the minimum error objective value achieved at each {bus location}, with the lowest value at bus 13 as the correct location.} \label{fig:obj_14}
\end{figure}

Upon solving the optimization problem \eqref{eq:mc_3} for every bus location, we plot the resultant minimum mismatch error in Fig. \ref{fig:obj_14}. As any {\yuqi bus split} at bus 8 leads to system islanding, it is excluded from the comparison. 
The mismatch error at bus 13 is the smallest among all possible buses, and thus it is identified as the correct location. Moreover, the optimal solution for this bus has $\bm{z} = [0;0;1]$, which correctly indicates that line $(13,14)$ is reconnected to the new bus while $z^d = 1$  identifies the reconnected load.
\begin{figure}[t!]
\centering
\vspace{-2pt}
\includegraphics[trim=0cm 0cm 1cm 0cm,clip=true,totalheight=0.16\textheight]{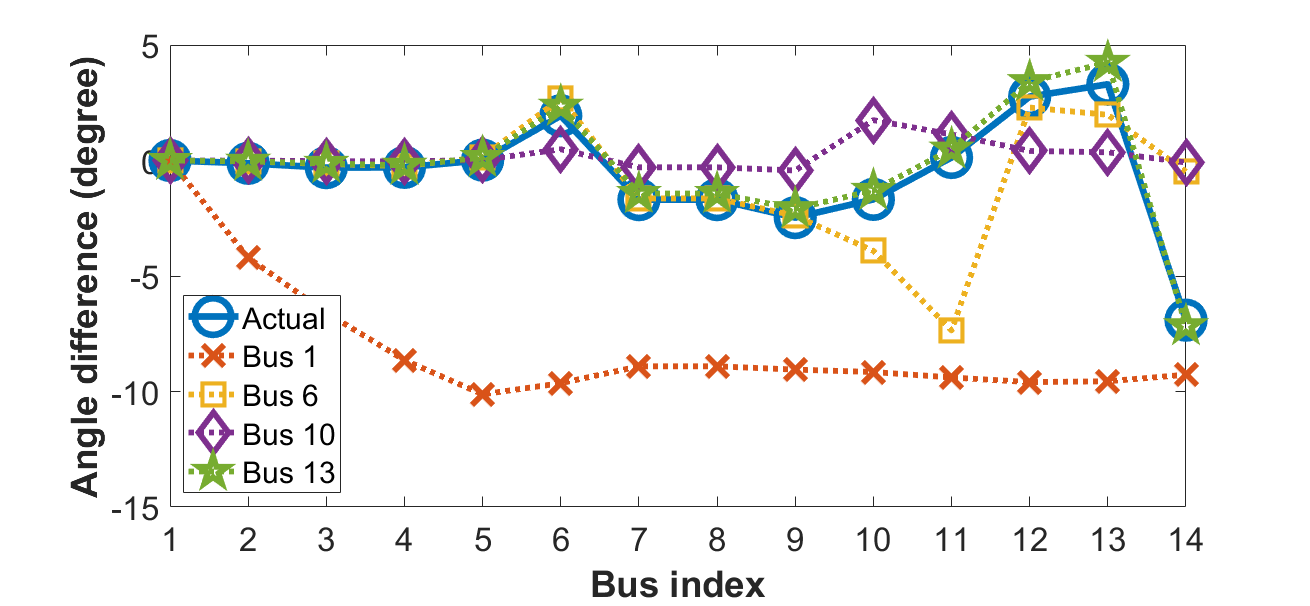}
\caption{System-wide phase angle difference given by selected bus locations as compared to the actual values.}
\label{fig:angle_comp_2}
\end{figure}
Using the optimal solution for each bus location, we can  find the corresponding system-wide phase angle difference based on dc model, namely 
$[\tilde{\mathbf{B}}(\bm{z})]^{-1} \tilde{\mathbf{p}}(\bm{z}) - \hat{\bm{\theta}}$.
Fig.\ \ref{fig:angle_comp_2} plots the resultant angle difference solutions for selected bus locations as compared to the actual values $\boldsymbol \delta$  from the ac power flow model as observed by PMUs.
This comparison again confirms that the split at bus 13 has the closest match with the actual system responses. 

It is worth mentioning that this test also points out how to efficiently select candidate bus locations for larger-sized systems. As shown by Fig. \ref{fig:angle_comp_2}, the actual phase angles at bus 13 and its neighbor buses $\mathcal{N}_{13} = \{6, 12, 14\}$ change more evidently than other buses. As the inverse of Bbus matrix is typically diagonally dominant, the effects of topology changes tend to reduce from the change location to buses further away; see e.g., \cite{zhu2014blocking}. This localized impact property enables to adopt a simple yet effective approach to prioritize candidate bus locations based on their angle differences, which will turn out to be very useful for large systems as detailed soon.

\subsection{IEEE 300-Bus System Tests}
\label{sec:num_300}
We have performed more comprehensive validations using the IEEE 300-bus system that provide quantifiable identification performance. We first investigate the bus split at the substation of bus $i=120$, and then provide more quantitative results for different bus locations. Finally, the impact of partial observation is also studied.

\textit{Test Case 1:} We pick the split of bus $i=120$ to investigate. It has 3 neighboring buses, and is connected with both generation and load. Compared to the 14-bus system, there are a larger number of buses in this system, making it difficult to enumerate all the buses as candidate buses. Using the empirical insights earlier on, we use a heuristic selection scheme that ranks the buses based on the observed value of phase angle difference, namely $|\delta_i|$ provided by the synchrophasor data.  
If not all buses are equipped with PMUs, as later in Test Case 3,
the neighboring buses of the highly-ranked buses would also be included as candidate locations. For the split of bus $i=120$, the top six buses in the ranking are buses 120, 153, 151, 152, 155, and 154, most of which are co-located in the same area to bus 120. These six buses are selected to run the problem \eqref{eq:mc_3}. Accordingly, the bus with the smallest achievable error objective value is deemed as the location of bus split, along with its corresponding optimal solution. 

\begin{figure}[t!]
\centering
\vspace{-2pt}
\includegraphics[trim=0cm 0cm 1cm 0cm,clip=true,totalheight=0.16\textheight]{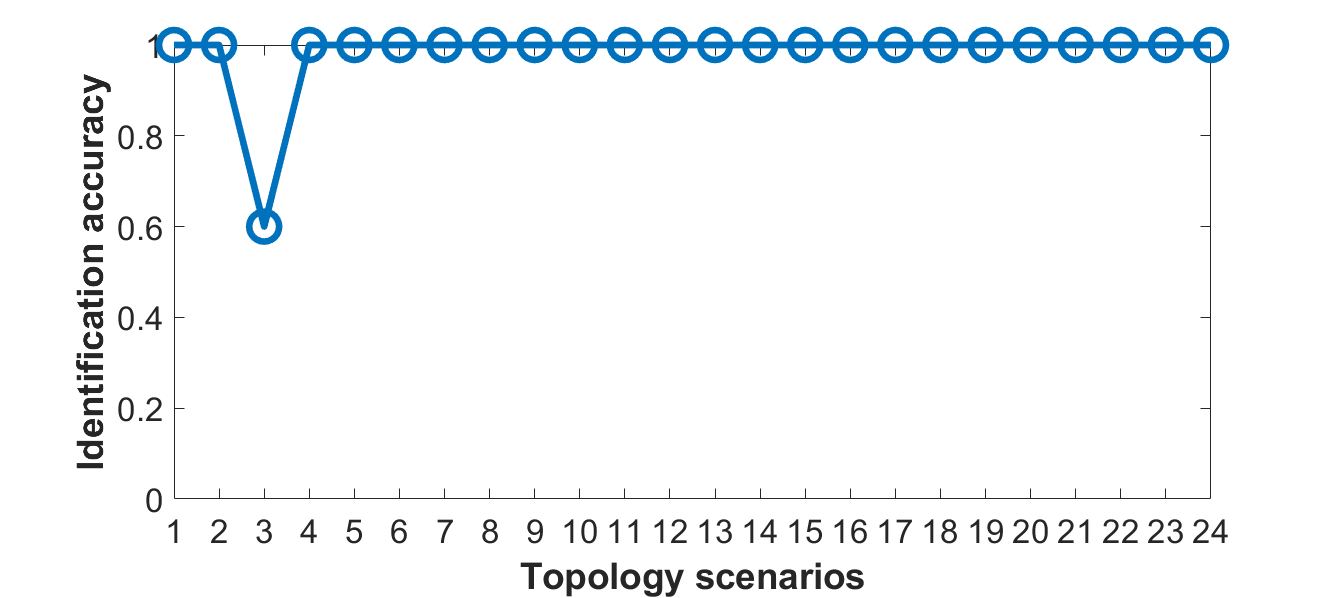}
\caption{Identification accuracy of different scenarios for the bus split in Test Case 1.}
\label{fig:accuracy_120}
\end{figure}

This identification process has been run for every possible topology scenario (every choice of binary $\bm{z}$) for the split of bus 120. The identification accuracy is evaluated based on the percentage of correctly recovered entries of vector $\bm{z}$, which is given in  Fig. \ref{fig:accuracy_120} for every topology scenario.
For a majority of scenarios, the identification accuracy is perfect. This implies that the angle differences for each of these scenarios sufficiently differentiate from those of all other scenarios, and our proposed solution can effectively find the correct scenario. Nonetheless, this is not the case for scenario 3, where the connections for two out of the three lines have been erroneously identified. 
A closer look at this scenario reveals that the dc approximation error has led to that the mis-identified $\bm{z}$ solution has the smallest mismatch error with the nonlinear ac model.
Overall, the proposed method can correctly find the optimal solution to \eqref{eq:enumeration_1}, while the accuracy of latter may still depend on the approximation error of the dc model.

\begin{figure}[t!]
\centering
\vspace{-2pt}
\includegraphics[trim=0cm 0cm 1.5cm 0cm,clip=true,totalheight=0.16\textheight]{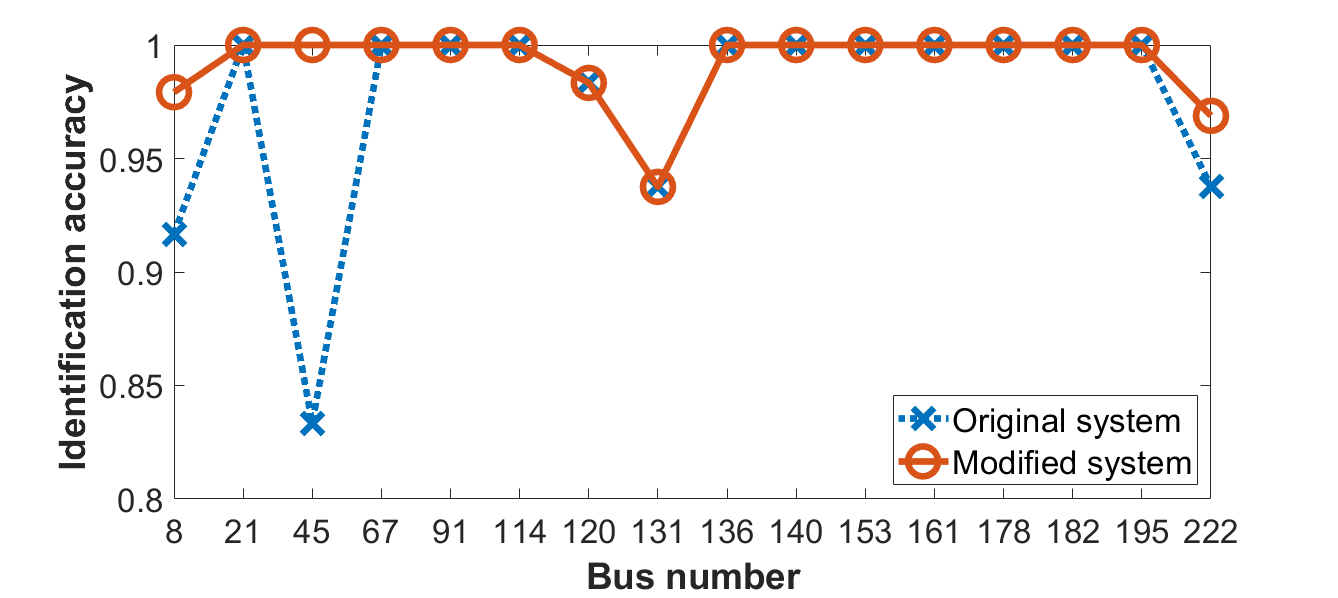}
\caption{Identification accuracy of bus split at selected {bus locations} in Test Case 2.}
\label{fig:accuracy_multiple}
\end{figure}

\textit{Test Case 2.} We have further tested our proposed identification algorithm for the bus split events at 16 selected buses, giving rise to 116 topology scenarios in total. To quantify the identification accuracy for each bus, we average the percentage of correct recovery over all the scenarios under the split of that bus. The resultant accuracy results for the original 300-bus system are given in Fig.  \ref{fig:accuracy_multiple}, along with those for a modified 300-bus system. Test Case 1 pointed out the impact of the dc power flow approximation, which is highly related to  the line resistance-reactance ratio \cite{stott2009dc}. Some of the transmission lines in the original 300-bus system have a quite high ratio value, which is uncommon for high-voltage grids. Hence, we reduce the resistance values therein by half to attain a modified system.
Overall, our proposed algorithm has achieved effective identification of bus split events in all the 16 selected buses. On average, the original system gives an accuracy of $97.6\%$, with a minimum accuracy of $83.3\%$ at bus 45. By decreasing the resistance values, the modified system enjoys an increased average accuracy of $99.2\%$. 
Hence, our identification algorithm can attain accurate results for practical systems. 

\begin{figure}[t!]
\centering
\vspace{-2pt}
\includegraphics[trim=0cm 0cm 1.5cm 0cm,clip=true,totalheight=0.16\textheight]{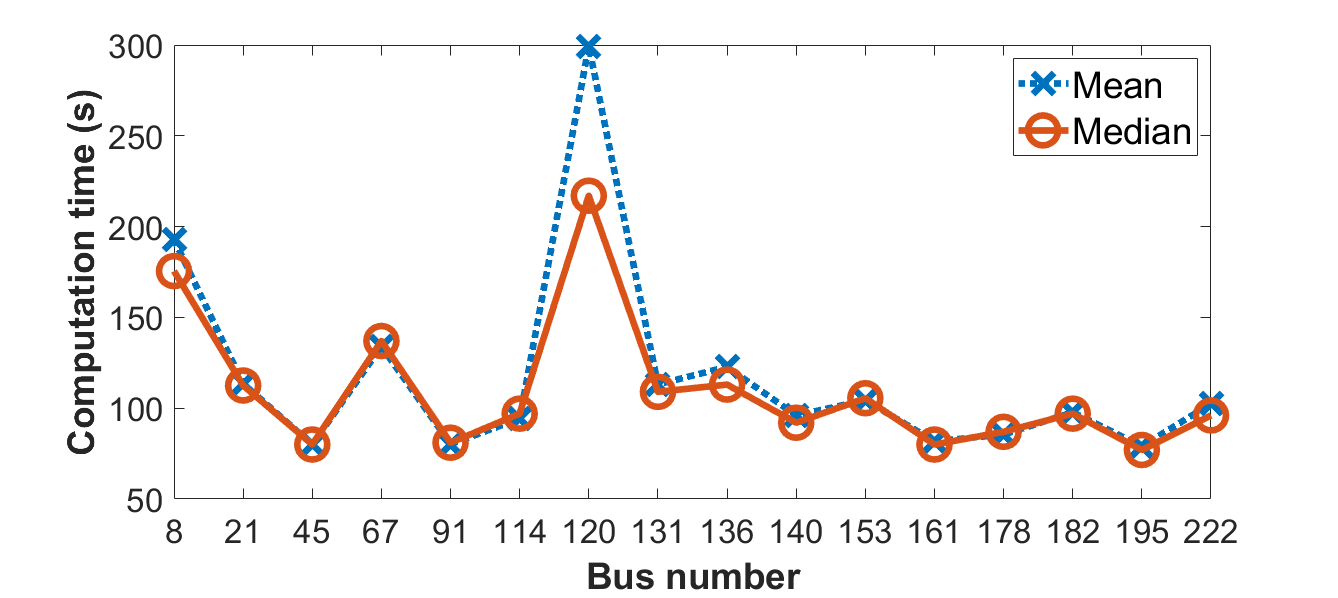}
\caption{Computation time for identification for selected {bus split locations} in Test Case 2.}
\label{fig:time}
\end{figure}

To evaluate the computational complexity of the proposed algorithm, the mean and median values for the run-time of all topology scenarios per bus for the original 300-bus system are given in Fig. \ref{fig:time}.
For this 300-bus system, it takes around 100 seconds to execute the proposed identification algorithm using a standard computer. This run-time scales nicely with the size of system. 
{Improved parameter settings (e.g., the upper bound ${X}^{\max}_{mn}$) and high-performance computing resources can further facilitate the implementation of the proposed algorithm for real-time monitoring of bus split events.}


\textit{Test Case 3.} This test investigates the effects of partial observation using the modified system of reduced resistance. As mentioned in Sec. \ref{sec:id}, the problem \eqref{eq:enumeration_1} can be easily modified to allow for partial angle or additional line flow measurements. Using the benchmark of full angle measurements, we compare three partial observation settings: $70\%$ angle measurements, $85\%$ angle measurements and a combination of $70\%$ angle measurements and $50\%$ line flow measurements. For the partial angle measurements, it suffices to use a selection matrix to extract the mismatch error vector in \eqref{eq:enumeration_1} for the metered locations. As for the line flow measurements, they can be incorporated into the error mismatch objective similar to angle measurements by recognizing the linear relation in \eqref{eq:line_flow}. 
All the measurement locations have been randomly selected and the mismatch error vector is normalized to match the scaling difference between angle and flow changes. Specifically, the ratio between the average of the absolute angle changes and that of the absolute line flow change is used to scale these two types of measurements.

\begin{figure}[t!]
\centering
\vspace{-2pt}
\includegraphics[trim=0cm 0cm 1.5cm 0cm,clip=true,totalheight=0.16\textheight]{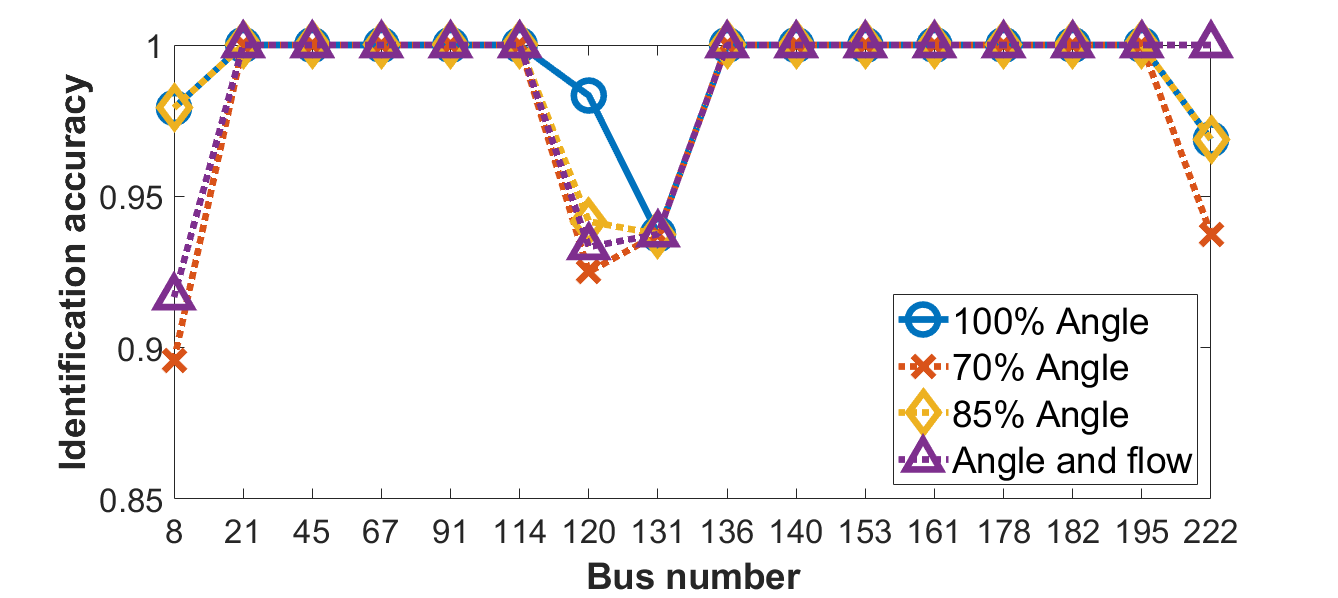}
\caption{Identification accuracy for three partial observation cases compared to the benchmark in Test Case 3.}
\label{fig:partial}
\end{figure}

Fig. \ref{fig:partial} plots the identification accuracy for the partial observation settings along with the benchmark. All
the three settings have achieved satisfactory accuracy (on average $98.1\%$, $98.9\%$ and $98.7\%$ for each of the three settings), with a minimum accuracy of $89.6\%$ at bus 8 under the $70\%$ angle setting. 
Generally speaking, the accuracy consistently increases when there are more measurements of angles or line flows. Note that for the split events at bus 222, the inclusion of line flow data is very useful for achieving perfect identification results. At some other buses (bus 8 or 120), the additional angle information is more helpful for identifying certain split events. 
Therefore, the proposed algorithm can adapt to various measurement availability conditions and provide an accurate identification performance for large systems.


\section{Conclusions} \label{sec:con}

This paper developed an efficient modeling and monitoring framework for power system bus split events due to substation connectivity changes. Based on the dc power flow analysis for the compact bus-branch model, the linear sensitivity analysis was performed that can quickly compute the system-wide changes for given bus split event. In addition, synchrophasor data enabled monitoring of bus split events was cast as an optimization problem, with binary variables indicating the connectivity of lines, generation, loads within a substation. To tackle the bilinear relations in the resultant problem, 
the McCormick relaxation technique has been leveraged to attain an \textit{equivalent} MILP reformulation that is efficiently solvable for real-time identification. Numerical studies have corroborated the performance of the proposed identification algorithm for enhancing power system situational awareness in the face of bus split events.

%

\bibliography{bibliography,hzpub}

\begin{thebibliography}{10}
\providecommand{\url}[1]{#1}
\csname url@samestyle\endcsname
\providecommand{\newblock}{\relax}
\providecommand{\bibinfo}[2]{#2}
\providecommand{\BIBentrySTDinterwordspacing}{\spaceskip=0pt\relax}
\providecommand{\BIBentryALTinterwordstretchfactor}{4}
\providecommand{\BIBentryALTinterwordspacing}{\spaceskip=\fontdimen2\font plus
\BIBentryALTinterwordstretchfactor\fontdimen3\font minus
  \fontdimen4\font\relax}
\providecommand{\BIBforeignlanguage}[2]{{%
\expandafter\ifx\csname l@#1\endcsname\relax
\typeout{** WARNING: IEEEtran.bst: No hyphenation pattern has been}%
\typeout{** loaded for the language `#1'. Using the pattern for}%
\typeout{** the default language instead.}%
\else
\language=\csname l@#1\endcsname
\fi
#2}}
\providecommand{\BIBdecl}{\relax}
\BIBdecl

\bibitem{wood2013power}
A.~J. Wood, B.~F. Wollenberg, and G.~B. Shebl{\'e}, \emph{Power generation,
  operation, and control}.\hskip 1em plus 0.5em minus 0.4em\relax John Wiley \&
  Sons, 2013.

\bibitem{EPRI01}
S.~Santoso, ``Transmission line protection for enhanced security: Including
  circuit models in relays,'' EPRI, Tech. Rep., 2014.

\bibitem{case2016analysis}
R.~M. Lee, M.~J. Assante, and T.~Conway, ``Analysis of the cyber attack on the
  {U}krainian power grid,'' Electricity Information Sharing and Analysis Center
  (E-ISAC), Tech. Rep., 2016.

\bibitem{abur1995identifying}
A.~Abur, H.~Kim, and M.~Celik, ``Identifying the unknown circuit breaker
  statuses in power networks,'' \emph{IEEE Transactions on Power Systems},
  vol.~10, no.~4, pp. 2029--2037, 1995.

\bibitem{kassakian2011future}
J.~G. Kassakian, R.~Schmalensee, G.~Desgroseilliers \emph{et~al.}, ``The future
  of the electric grid,'' Massachusetts Institute of Technology, Tech. Rep.,
  2011.

\bibitem{WECC1}
{Contingency Subgroup of the Modeling SPS and Relays Ad-Hoc Task Force},
  ``Node-breaker white paper,'' WECC, Tech. Rep., 2014.

\bibitem{NERC2}
``Node-breaker modeling representation,'' NERC, Tech. Rep., 2016.

\bibitem{NERC1}
{NERC Planning Committees}, ``Proposal for development and use of node breaker
  topology representations for off-line study models,'' NERC, Tech. Rep., 2019.

\bibitem{abur2004power}
A.~Abur and A.~G. Exposito, \emph{Power system state estimation: theory and
  implementation}.\hskip 1em plus 0.5em minus 0.4em\relax CRC press, 2004.

\bibitem{wu1989detection}
F.~F. Wu and W.-H. Liu, ``Detection of topology errors by state estimation
  (power systems),'' \emph{IEEE Transactions on Power Systems}, vol.~4, no.~1,
  pp. 176--183, 1989.

\bibitem{clements1998topology}
K.~A. Clements and A.~S. Costa, ``Topology error identification using
  normalized lagrange multipliers,'' \emph{IEEE Transactions on Power Systems},
  vol.~13, no.~2, pp. 347--353, 1998.

\bibitem{exposito2001reduced}
A.~G. Exp{\'o}sito and A.~de~la Villa~Jaen, ``Reduced substation models for
  generalized state estimation,'' \emph{IEEE Transactions on Power Systems},
  vol.~16, no.~4, pp. 839--846, 2001.

\bibitem{korres2002identification}
G.~N. Korres and P.~J. Katsikas, ``Identification of circuit breaker statuses
  in {WLS} state estimator,'' \emph{IEEE Transactions on Power Systems},
  vol.~17, no.~3, pp. 818--825, 2002.

\bibitem{lourencco2004bayesian}
E.~M. Louren{\c{c}}o, A.~S. Costa, and K.~A. Clements, ``Bayesian-based
  hypothesis testing for topology error identification in generalized state
  estimation,'' \emph{IEEE Transactions on Power Systems}, vol.~19, no.~2, pp.
  1206--1215, 2004.

\bibitem{lourencco2014topology}
E.~M. Louren{\c{c}}o, E.~P. Coelho, and B.~C. Pal, ``Topology error and bad
  data processing in generalized state estimation,'' \emph{IEEE Transactions on
  Power Systems}, vol.~30, no.~6, pp. 3190--3200, 2014.

\bibitem{tate2008line}
J.~E. Tate and T.~J. Overbye, ``Line outage detection using phasor angle
  measurements,'' \emph{IEEE Transactions on Power Systems}, vol.~23, no.~4,
  pp. 1644--1652, 2008.

\bibitem{zhu2012sparse}
H.~Zhu and G.~B. Giannakis, ``Sparse overcomplete representations for efficient
  identification of power line outages,'' \emph{IEEE Transactions on Power
  Systems}, vol.~27, no.~4, pp. 2215--2224, 2012.

\bibitem{emami2012external}
R.~Emami and A.~Abur, ``External system line outage identification using phasor
  measurement units,'' \emph{IEEE Transactions on Power Systems}, vol.~28,
  no.~2, pp. 1035--1040, 2012.

\bibitem{chen2015quickest}
Y.~C. Chen, T.~Banerjee, A.~D. Dominguez-Garcia, and V.~V. Veeravalli,
  ``Quickest line outage detection and identification,'' \emph{IEEE
  Transactions on Power Systems}, vol.~31, no.~1, pp. 749--758, 2015.

\bibitem{deka2016estimating}
D.~Deka, S.~Backhaus, and M.~Chertkov, ``Estimating distribution grid
  topologies: A graphical learning based approach,'' in \emph{2016 Power
  Systems Computation Conference (PSCC)}.\hskip 1em plus 0.5em minus
  0.4em\relax IEEE, 2016, pp. 1--7.

\bibitem{cavraro2018graph}
G.~Cavraro and V.~Kekatos, ``Graph algorithms for topology identification using
  power grid probing,'' \emph{IEEE control systems letters}, vol.~2, no.~4, pp.
  689--694, 2018.

\bibitem{ardakanian2019identification}
O.~Ardakanian, V.~W. Wong, R.~Dobbe, S.~H. Low, A.~von Meier, C.~J. Tomlin, and
  Y.~Yuan, ``On identification of distribution grids,'' \emph{IEEE Transactions
  on Control of Network Systems}, vol.~6, no.~3, pp. 950--960, 2019.

\bibitem{wu2016online}
M.~Wu and L.~Xie, ``Online detection of low-quality synchrophasor measurements:
  A data-driven approach,'' \emph{IEEE Transactions on Power Systems}, vol.~32,
  no.~4, pp. 2817--2827, 2016.

\bibitem{zhou2018ensemble}
M.~Zhou, Y.~Wang, A.~K. Srivastava, Y.~Wu, and P.~Banerjee, ``Ensemble-based
  algorithm for synchrophasor data anomaly detection,'' \emph{IEEE Transactions
  on Smart Grid}, vol.~10, no.~3, pp. 2979--2988, 2018.

\bibitem{ramakrishna2019detection}
R.~Ramakrishna and A.~Scaglione, ``Detection of false data injection attack
  using graph signal processing for the power grid,'' in \emph{2019 IEEE Global
  Conference on Signal and Information Processing (GlobalSIP)}.\hskip 1em plus
  0.5em minus 0.4em\relax IEEE, 2019, pp. 1--5.

\bibitem{wang2019data}
S.~Wang, P.~Dehghanian, and B.~Zhang, ``A data-driven algorithm for online
  power grid topology change identification with {PMU}s,'' in \emph{Proc. IEEE
  PES General Meeting}, 2019.

\bibitem{heidarifar2015network}
M.~Heidarifar and H.~Ghasemi, ``A network topology optimization model based on
  substation and node-breaker modeling,'' \emph{IEEE Transactions on Power
  Systems}, vol.~31, no.~1, pp. 247--255, 2015.

\bibitem{park2019sparse}
B.~Park, J.~Holzer, and C.~L. DeMarco, ``A sparse tableau formulation for
  node-breaker representations in security-constrained optimal power flow,''
  \emph{IEEE Transactions on Power Systems}, vol.~34, no.~1, pp. 637--647,
  2019.

\bibitem{donmez2020parallel}
B.~Donmez and A.~Abur, ``A parallel framework for robust state estimation using
  node-breaker substation models,'' in \emph{2020 IEEE PES Innovative Smart
  Grid Technologies Europe (ISGT-Europe)}.\hskip 1em plus 0.5em minus
  0.4em\relax IEEE, 2020, pp. 1136--1140.

\bibitem{kekatos2012joint}
V.~Kekatos and G.~B. Giannakis, ``Joint power system state estimation and
  breaker status identification,'' in \emph{Proc. North American Power Symp.},
  2012.

\bibitem{stott2009dc}
B.~Stott, J.~Jardim, and O.~Alsa{\c{c}}, ``{DC} power flow revisited,''
  \emph{IEEE Transactions on Power Systems}, vol.~24, no.~3, pp. 1290--1300,
  2009.

\bibitem{McCormick1976computability}
G.~P. McCormick, ``Computability of global solutions to factorable nonconvex
  programs: {Part I-Convex underestimating problems},'' \emph{Mathematical
  programming}, vol.~10, no.~1, pp. 147--175, 1976.

\bibitem{stott1974fast}
B.~Stott and O.~Alsac, ``Fast decoupled load flow,'' \emph{IEEE Transactions on
  Power Apparatus and Systems}, no.~3, pp. 859--869, 1974.

\bibitem{baldick2003variation}
R.~Baldick, ``Variation of distribution factors with loading,'' \emph{IEEE
  Transactions on Power Systems}, vol.~18, no.~4, pp. 1316--1323, 2003.

\bibitem{santoso2018standard}
S.~Santoso and H.~W. Beaty, \emph{Standard Handbook for Electrical
  Engineers}.\hskip 1em plus 0.5em minus 0.4em\relax McGraw-Hill Education,
  2018.

\bibitem{deka2015one}
D.~Deka, R.~Baldick, and S.~Vishwanath, ``One breaker is enough: Hidden
  topology attacks on power grids,'' in \emph{Proc. IEEE PES General Meeting},
  2015.

\bibitem{zhou2018false}
Y.~Zhou, J.~Cisneros-Saldana, and L.~Xie, ``False analog data injection attack
  towards topology errors: Formulation and feasibility analysis,'' in
  \emph{Proc. IEEE PES General Meeting}, 2018.

\bibitem{sun2003splitting}
K.~Sun, D.-Z. Zheng, and Q.~Lu, ``Splitting strategies for islanding operation
  of large-scale power systems using {OBDD}-based methods,'' \emph{IEEE
  Transactions on Power Systems}, vol.~18, no.~2, pp. 912--923, 2003.

\bibitem{kyriacou2017controlled}
A.~Kyriacou, P.~Demetriou, C.~Panayiotou, and E.~Kyriakides, ``Controlled
  islanding solution for large-scale power systems,'' \emph{IEEE Transactions
  on Power Systems}, vol.~33, no.~2, pp. 1591--1602, 2017.

\bibitem{boyd2004convex}
S.~Boyd and L.~Vandenberghe, \emph{Convex optimization}.\hskip 1em plus 0.5em
  minus 0.4em\relax Cambridge, U.K.: Cambridge University Press, 2004.

\bibitem{zhou2019bus}
Y.~Zhou and H.~Zhu, ``Bus split sensitivity analysis for enhanced security in
  power system operations,'' in \emph{Proc. North American Power Symp.}, 2019.

\bibitem{sherman1950adjustment}
J.~Sherman and W.~J. Morrison, ``Adjustment of an inverse matrix corresponding
  to a change in one element of a given matrix,'' \emph{The Annals of
  Mathematical Statistics}, vol.~21, no.~1, pp. 124--127, 1950.

\bibitem{dorfler2012kron}
F.~Dorfler and F.~Bullo, ``Kron reduction of graphs with applications to
  electrical networks,'' \emph{IEEE Transactions on Circuits and Systems I:
  Regular Papers}, vol.~60, no.~1, pp. 150--163, 2012.

\bibitem{gupte2013solving}
A.~Gupte, S.~Ahmed, M.~S. Cheon, and S.~Dey, ``Solving mixed integer bilinear
  problems using {MILP} formulations,'' \emph{SIAM Journal on Optimization},
  vol.~23, no.~2, pp. 721--744, 2013.

\bibitem{kocuk2017new}
B.~Kocuk, S.~S. Dey, and X.~A. Sun, ``New formulation and strong {MISOCP}
  relaxations for {AC} optimal transmission switching problem,'' \emph{IEEE
  Transactions on Power Systems}, vol.~32, no.~6, pp. 4161--4170, 2017.

\bibitem{bhela2019designing}
S.~Bhela, D.~Deka, H.~Nagarajan, and V.~Kekatos, ``Designing power grid
  topologies for minimizing network disturbances: An exact {MILP}
  formulation,'' in \emph{Proc. American Control Conference (ACC)}, 2019.

\bibitem{park2020optimal}
B.~Park and C.~L. Demarco, ``Optimal network topology for node-breaker
  representations with ac power flow constraints,'' \emph{IEEE Access}, vol.~8,
  pp. 64\,347--64\,355, 2020.

\bibitem{ostrowski2011tight}
J.~Ostrowski, M.~F. Anjos, and A.~Vannelli, ``Tight mixed integer linear
  programming formulations for the unit commitment problem,'' \emph{IEEE
  Transactions on Power Systems}, vol.~27, no.~1, pp. 39--46, 2011.

\bibitem{plemmons1977m}
R.~J. Plemmons, ``M-matrix characterizations. {I}-nonsingular m-matrices,''
  \emph{Linear Algebra and its Applications}, vol.~18, no.~2, pp. 175--188,
  1977.

\bibitem{zhu2014blocking}
H.~Zhu and T.~J. Overbye, ``Blocking device placement for mitigating the
  effects of geomagnetically induced currents,'' \emph{IEEE Transactions on
  Power Systems}, vol.~30, no.~4, pp. 2081--2089, 2014.

\end{thebibliography}





\bibliographystyle{IEEEtran}

\itemsep2pt

\end{document}